\documentclass[12pt]{amsart}
\usepackage[myheadings]{fullpage}
\usepackage{amsmath,amsthm}                                                     

\title{Volume comparison via boundary distances}
\author{Sergei Ivanov}
\thanks{Supported by the Dynasty foundation and RFBR grants
08-01-00079-a, 09-01-12130-ofi-m}
\email{svivanov@pdmi.ras.ru}
\address{St.~Petersburg Department of Steklov Institute of Mathematics,
Fontanka 27, 191023, St.~Petersburg, Russia}

\numberwithin{equation}{section}

\newtheorem{theorem}{Theorem}[section]
\newtheorem{conjecture}[theorem]{Conjecture}
\newtheorem{prop}[theorem]{Proposition}

\theoremstyle{definition}
\newtheorem{definition}[theorem]{Definition}
\newtheorem{question}[theorem]{Question}

\newtheorem*{remark}{Remark}

\newcommand{\R}{\mathbb R}
\newcommand{\Z}{\mathbb Z}
\newcommand{\HH}{\mathbb H}

\DeclareMathOperator{\FillVol}{FillVol}
\DeclareMathOperator{\AsVol}{AsVol}
\DeclareMathOperator{\vol}{Vol}
\DeclareMathOperator{\area}{area}

\begin{document}

\begin{abstract}
The main subject of this lecture is a connection
between Gromov's filling volumes
and a boundary rigidity problem of determining
a Riemannian metric in a compact domain
by its boundary distance function.
A fruitful approach is to represent Riemannian metrics
by minimal surfaces in a Banach space and to prove rigidity
by studying the equality case in a filling volume inequality.
I discuss recent results obtained with this approach
and related problems in Finsler geometry.
\end{abstract}

\maketitle

\section{Introduction}

\subsection{A toy question}

One of the goals of this lecture is to advertise
a conjecture about filling volumes. It can be
stated without preliminaries
(although in an obscured way) as follows.

\begin{question}\label{q-basic}
Let $N^{n+1}$ be a complete Riemannian manifold and $M^n\subset N$
a compact hypersurface with boundary.
Suppose that $M$ is convex in the following strong sense:
for every two points $x,y\in M$, there is a unique shortest
geodesic segment connecting $x$ and~$y$ in $N$, and this segment
lies in~$M$.
(In particular, $M$ is totally geodesic.)

Is it true that every such $M$ is an area minimizer?
That is, does it have the least $n$-dimensional
area among all compact (orientable)
hypersurfaces in $N$ with the same boundary?
\end{question}

The wording of this question is deliberately chosen so as
to make an affirmative answer sound more plausible.
Actually the answer is not known, and
an affirmative one would have strong implications.

The convexity assumptions in Question \ref{q-basic}
imply that $M$ is diffeomorphic to the $n$-disc, its boundary is
convex, and all its geodesic are shortest paths.
The latter is a crucial property while the former two could be relaxed:
for example, non-convex regions in $M$ are area minimizers if so is $M$.

Since the surface in question is totally geodesic,
it is minimal in the variational sense: the mean curvature,
and hence the first variation of area, is zero.
Cutting off a neighborhood of the boundary yields
a surface where geodesics have no conjugate points,
and it is easy to see that in this case it is
a stable minimal surface and hence minimizes the area
locally (among all nearby surfaces).
However the \textit{global} area-minimality
in Question~\ref{q-basic} is a completely different issue.

\subsection{Boundary rigidity}
I postpone further discussion of Question \ref{q-basic}
until subsection \ref{sec-filling}.
This subsection is a brief introduction to boundary distance rigidity.

For a Riemannian manifold $M$, possibly with boundary,
let $d_M$ denote the induced length metric on~$M$.
This is a function on $M\times M$ measuring geodesic
distances between points.
The \textit{boundary distance function} of $M$, denoted by $bd_M$,
is the restriction of $d_M$ to $\partial M\times\partial M$.
It is natural to ask whether the metric in the interior can be
determined if one knows the boundary distance function.

Inverse boundary problems of this type were originally
motivated by geophysics: the inner structure of the Earth
can be studied by measuring travel times of seismic waves
between points at the surface.
Assuming that the Earth is filled by
isotropic media with variable speed of sound,
the travel times represent the boundary distance function of a
conformal metric on~$D^3$, and the problem is to determine the conformal
factor by these data.
Under the assumption that the Earth is spherically symmetric,
this inverse kinematic problem was solved by Herglotz \cite{herglotz} and
Wiechert~\cite{wiechert}. For a general simple conformal metric, the uniqueness
of a solution was proved by Mukhometov and Romanov \cite{Mu-Ro},
see also \cite{beylkin}, \cite{Mu}, \cite{croke91}.

If the metric is not supposed to be conformal, determining metric coefficients
as functions of coordinates does not make sense: any Riemannian isometry that fixes
the boundary obviously preserves the boundary distances. Two metrics
related by such an isometry must be regarded as the same metric, hence
the following definition.

\begin{definition}\label{d-rigid}
A compact Riemannian manifold $M$ with boundary is said to be \textit{boundary rigid}
if it is determined by its boundary distance function uniquely up to an isometry
fixing the boundary.

In a more formal language this means the following:
every compact Riemannian manifold $M'$ such that
$\partial M'=\partial M$ and $bd_{M'}=bd_M$ is isometric to $M$ via an isometry
$f:M\to M'$ such that $f|_{\partial M} =id_{\partial M}$.
\end{definition}

It is easy to construct metrics that are \emph{not} boundary rigid.
For example, begin with an arbitrary metric and enlarge it
near a point $p$ so that no shortest path between boundary points
goes through~$p$. Then a perturbation of the metric
near $p$ does not affect the boundary distance function.
Another example is the standard hemisphere: since the boundary
distances are realized by boundary arcs, enlarging the metric
in the interior does not change them.

Such examples should be excluded if one seeks boundary rigidity.
A natural set of restrictions is contained in the following definition.

\begin{definition}\label{d-simple}
A compact Riemannian manifold $M$ is said to be \textit{simple} if

(1) The boundary $\partial M$ is strictly convex,
i.e.\ has positive definite second fundamental form.

(2) Every geodesic segment in $M$ is minimal, i.e.\ realizes
the distance between its endpoints.

(3) The geodesics in $M$ have no conjugate points.
(Or, equivalently, there is a larger manifold $M^+$ containing $M$
in its interior and such that all geodesics in $M^+$ are minimal.)
\end{definition}

For example, the standard hemisphere is not simple
but cutting off an arbitrarily small neighborhood
of the boundary makes it simple.

The first requirement of Definition \ref{d-simple}
implies that all distances in $M$
are realized by geodesics. Then one easily sees that the
exponential map at every point is a diffeomorphism,
and it follows that a simple manifold is diffeomorphic to
a disc. Thus one may as well speak about \textit{simple metrics} on~$D^n$.

Note that simplicity of the metric can be observed via the boundary
distance function. That is, if two metrics have the same boundary
distance function, then either they are both simple or both are not. 
Indeed, the convexity of $\partial M$ is equivalent to a sort of strict
triangle inequality for $bd_M$, and the fact that geodesics are minimal
and have no conjugate points is equivalent to smoothness of
$bd_M$ away from the diagonal.

\begin{conjecture}[R.~Michel \cite{Michel}]
\label{c-michel}
Every simple Riemannian manifold is boundary rigid.
\end{conjecture}

Pestov and Uhlmann \cite{pestov-uhlmann} proved this
conjecture in dimension~2. In higher dimensions
the following types of spaces are known to be
boundary rigid:
\begin{itemize}
 \item regions in $\R^n$ and moreover all $n$-dimensional
flat manifolds that admit an isometric immersion to $\R^n$
(Besikovitch \cite{Bes}; Gromov \cite{gromov-frm});
 \item regions in the standard open hemisphere $S^n_+$
(Michel \cite{Michel});
 \item regions is symmetric spaces of negative curvature
(this follows from a volume entropy inequality
proved by Besson, Courtois and Gallot \cite{BCG});
 \item regions in metric products of the form $M_0\times\R$
where $M_0$ is a complete simply connected
Riemannian manifold without conjugate points
(Croke and Kleiner \cite{CK98});
 \item metrics sufficiently close in $C^2$ to the Euclidean
metric of a region in $\R^n$ (Burago and Ivanov \cite{BI10});
 \item metrics sufficiently close in $C^3$
to the hyperbolic metric of a region in $\HH^n$
(Burago and Ivanov \cite{BI-new}).
\end{itemize}
Proofs of the last two results are discussed in section \ref{sec-banach}.

\begin{remark}
More is known about the local variant of the conjecture,
that is, when the metrics of $M$ and $M'$ in Definition \ref{d-rigid}
are assumed \textit{a priori} close to each other.
Local boundary rigidity is proved for a generic set of simple metrics
including all analytic ones \cite{SU} and for all metrics
with ``not too much'' positive curvature \cite{CDS}.
\end{remark}

\subsection{Filling volumes and minimal fillings}
\label{sec-filling}

To simplify matters, all manifolds and surfaces in the sequel are assumed
orientable.
And for the most part one may assume that all Riemannian manifolds
in question are just metrics on the disc~$D^n$.

\begin{definition}
\label{d-filvol}
Let $N$ be a closed $(n-1)$-dimensional manifold and $f:N\times N\to\R$
a nonnegative function. The \textit{filling volume} of $f$,
denoted by $\FillVol(N,f)$, is defined by
\begin{equation}
\label{e-fillvol}
 \FillVol(N,f) = \inf \{ \vol(M): \partial M=N,\ bd_M\ge f \}
\end{equation}
where the infimum is taken over all (orientable) compact
$n$-dimensional Riemannian manifolds $M$
such that $\partial M=N$ and $bd_M\ge f$.
Such manifolds $M$ are referred to as \textit{fillings}
of $(N,f)$.

A compact Riemannian manifold $M$ is said to be a \textit{minimal filling}
if it realizes the infimum in \eqref{e-fillvol} for $S=\partial M$
and some function $f$ (and hence for $f=bd_M$).
In other words, $M$ is a minimal filling if
$\vol(M) = \FillVol(\partial M,bd_M)$.
\end{definition}

The notion of filling volume was
introduced by Gromov \cite{gromov-frm}
(in the case when $f$ is a metric on~$N$).
The above definition assumes that there are no
topological obstructions for $N$ to be a boundary,
cf.~\cite{gromov-frm} for the general case.

Substituting intermediate definitions yields the following:
$M$ is a minimal filling if and only if, for every
compact Riemannian manifold $M'$ such that
$\partial M'=\partial M$ and
\begin{equation}
\label{e-fillboundary}
 d_{M'}(x,y)\ge d_M(x,y)
\quad\text{for all $x,y\in\partial M$} ,
\end{equation}
one has
\begin{equation}
\label{e-minfill}
 \vol(M') \ge \vol(M) .
\end{equation}

The following conjecture is the main topic of this lecture.

\begin{conjecture}
\label{c-main}
Every simple manifold is a minimal filling.
\end{conjecture}

Note that a $(C^0)$ limit of minimal fillings is also
a minimal filling, and a limit of simple metrics can
have a non-strictly convex boundary and non-strictly minimal
geodesics. Thus the simplicity assumption in
Conjecture \ref{c-main} can be relaxed to allow for
such cases. In particular, if the conjecture is true,
then the standard hemisphere is a minimal filling.

Convexity of the boundary is a convenience assumption
and it can be removed in some cases (see e.g.~\cite{I09-pre}).
Observe that
any subregion of a minimal filling is a minimal filling
as well.

If a simple manifold $M$ is found to be a minimal filling,
one can try to analyze the equality case in \eqref{e-minfill}
and hope that it is attained only if $M'$ is isometric to~$M$
(via an isometry fixing the boundary).
This hope is expressed in the following
stronger variant of Conjecture \ref{c-main}.

\medbreak\noindent{\bf Conjecture \ref{c-main}$^{\mathbf+}$.}
{\it
Every simple manifold is a unique minimal filling of
its boundary distance function, up to an isometry
fixing the boundary.
\par}\medbreak

It is easy to see that Conjecture \ref{c-main}$^+$ implies
Michel's boundary rigidity conjecture \ref{c-michel}.
Almost all boundary rigid metrics listed above
are also known to be minimal fillings
(the exceptions are subsets of the hemisphere
and product metrics).
In dimension~2, all simple manifolds
are minimal fillings within the class of manifolds
homeomorphic to the disc \cite{I01}, but the general
filling minimality is not known even for the hemisphere.

Conjecture \ref{c-main} is equivalent to the affirmative
answer to Question \ref{q-basic}. Indeed, let $M\subset N$
be as in Question \ref{q-basic} and suppose that there
is a surface $M'\subset N$ with the same boundary but
smaller area. Then $M$ and $M'$, regarded as Riemannian
$n$-manifolds, satisfy \eqref{e-fillboundary} and hence
provide a counterexample to Conjecture \ref{c-main}.
Conversely, if manifolds $M$ and $M'$ satisfy
\eqref{e-fillboundary} but do not satisfy \eqref{e-minfill},
one can glue them together along the boundary and embed
the resulting space into a suitable manifold $N^{n+1}$
in order to produce a counterexample to Question \ref{q-basic}.
(One may need to change the metric of $M'$ near the boundary
to make a smooth gluing; this and other technical details
are easy to handle.)

\section{Some implications}

In this section I discuss some
implications of the minimal filling conjectures.

\subsection{Boundary rigidity}
As I already mentioned, Conjecture \ref{c-main}$^+$
implies Conjecture \ref{c-michel}.
Moreover, this implication works for every
individual manifold:

\begin{prop}
If a simple Riemannian manifold $M$ is a unique
minimal filling of its boundary distance function,
then $M$ is boundary rigid.
\end{prop}

The key to the proof is Santal\'o's integral geometric formula
for the volume of a simple Riemannian
manifold in terms of its boundary distance function and
its first order derivatives (cf.\ \cite{Santalo}, \cite{gromov-frm}, \cite{croke91}).
This formula implies that two simple manifolds with the same
boundary distance function have the same volume.
Recall that if $M$ is simple and $M'$ has the same
boundary distance function, then $M'$ is simple as well,
hence $\vol(M')=\vol(M)$ by Santal\'o's formula.
Then the uniqueness assumption
implies that $M$ and $M'$ are isometric.

This argument actually works not only for simple manifolds
but for a large class of \textit{strong geodesic minimizing} (SGM)
manifolds, cf.~\cite{croke91}.

\subsection{Gromov's circle filling conjecture}
What is the filling volume of the intrinsic metric of the circle?
This was the first question asked by Gromov
after the definition of filling volume in~\cite{gromov-frm}.
It is conjectured that this filling volume equals $2\pi$,
the value realised by the standard round hemisphere.
In other words, the question is: is the hemisphere
a minimal filling?
Since the hemisphere is a limit of simple manifolds,
Conjecture \ref{c-main} would immediately imply the
affirmative answer.

With definitions substituted,
the circle filling conjecture boils down to the following.
Let $M$ be a compact orientable two-dimensional surface with
a Riemannian metric such that $\partial M$ is a circle of length $2\pi$,
and for every pair $x,y$ of opposite points of this circle
one has ${d_{M}(x,y)=\pi}$. Then (the conjecture asserts that)
$\area(M)\ge 2\pi$.

This inequality is well-known if $M$ is homeomorphic to~$D^2$.
In other words, the hemisphere is a minimal
filling \textit{within the class of surfaces homeomorphic to the disc}.
Indeed, one can identify opposite points of the boundary circle
and obtain a closed surface $M_1\simeq\mathbb{RP}^2$ such that
the length of a shortest non-contractible loop in $M_1$ equals~$\pi$.
Then Pu's isosystolic inequality \cite{Pu} implies that
$\area(M)=\area(M_1)\ge 2\pi$.

Pu's original proof uses uniformization and integral geometry,
another proof can be found in \cite{I01}.
The uniformization approach can be pushed further
to cover the case when $M$ has genus~1,
cf.~\cite{BCIK05}. The case of a higher genus remains open.

The general case of the circle filling conjecture can
be similarly reformulated in terms of a systolic inequality,
and it has applications in higher-dimensional systolic geometry,
see e.g.~\cite[\S8.3]{katz}.

\subsection{E.~Hopf's theorem}
If $M$ is an $n$-torus with a Riemannian metric without conjugate
points, then $M$ is flat (that is, locally isometric to $\R^n$).
This fact was proved for $n=2$ by E.~Hopf \cite{hopf} and
for all $n$ by Burago and Ivanov~\cite{BI94}.
Both proofs involve dynamical arguments.
Croke and Kleiner \cite{CK95} proposed a more geometric approach
where E.~Hopf's theorem is derived from asymptotic volume
inequalities. Their approach led to a new proof of the theorem
in the two-dimensional case.
The following modification of their argument shows
how the theorem (in all dimensions) follows from Conjecture~\ref{c-main}
(with a relaxed boundary convexity assumption).

Let $\widetilde M$ denote the universal cover of $M$ with
the metric lifted from~$M$. The \textit{asymptotic volume} of $M$
is defined by
$$
 \AsVol(M) = \liminf_{R\to\infty} \frac{\vol(B_R)}{R^n}
$$
where $B_R$ is the metric ball in $\widetilde M$ centered at
a fixed point $x_0\in\widetilde M$.
Let $\omega_n$ denote the Euclidean volume of a unit ball in $\R^n$.
It can be shown that
\begin{equation}
\label{e-asvol}
\AsVol(M)\ge\omega_n ,
\end{equation}
with equality if and only if $M$ is flat.
This is proved in \cite{C92} for any closed Riemannian
$n$-manifold without conjugate points and
in \cite{BI95} for a Riemannian $n$-torus (with or without conjugate points).

Actually the inequality \eqref{e-asvol} can be improved by
inserting a factor depending on the affine type of the
stable norm $\|\cdot\|$ of~$M$, see \cite{BI95} and \cite[pp.~259--260]{gromov-metric}.
Namely
\begin{equation}
\label{e-asvol1}
 \AsVol(M)\ge \frac{\vol(B)}{\vol(E)}\cdot\omega_n
\end{equation}
where $B$ is the unit ball of $\|\cdot\|$ and $E$ is
the ellipsoid of maximal volume contained in $B$.
The equality in \eqref{e-asvol1} is attained if and only if
the metric is flat.

The universal cover $\widetilde M$ can be identified with $\R^n$
equipped with a $\Z^n$-periodic Riemannian metric.
Then the distances in $\widetilde M$ differ from the
distances in the normed space $(\R^n,\|\cdot\|)$ by
a bounded function, cf.~\cite{Bu92}. Let $d_E$ denote
the distance in the Euclidean metric associated with~$E$,
then 
\begin{equation}
\label{e-asvol2}
d_E(x,y)\ge d_{\widetilde M}(x,y)-\mbox{\it const}
\end{equation}
for all $x,y\in\widetilde M$. 
If $\widetilde M$ has no conjugate points, Conjecture \ref{c-main}
(without the boundary convexity assumption)
would imply that the ball $B_R\subset\widetilde M$ is a minimal filling.
Apply the minimal filling inequality \eqref{e-minfill} to $M=B_R$
and $M'=(B_R,d_E)$, where $d_E$ is modified near the boundary to
get rid of the constant in \eqref{e-asvol2}.
This yields the inequality opposite to \eqref{e-asvol1},
hence the metric of $\widetilde M$ is flat.

\section{Minimality in a Banach space}
\label{sec-banach}

In this section I discuss one of the approaches
to filling minimality and boundary rigidity
and outline the proofs of the following two theorems.

\begin{theorem}[\cite{BI10}]
\label{theorem1}
Let $D\subset\R^n$ be a compact region with a smooth boundary
and $g_0$ the standard Euclidean metric on $D$.
Then there is a neighborhood $\mathcal U$ of $g_0$
in the space of Riemannian metrics on~$D$
such that for every metric $g\in\mathcal U$
the space $(D,g)$ is a minimal filling and boundary rigid.
\end{theorem}

\begin{theorem}[\cite{BI-new}]
\label{theorem2}
Let $D\subset\HH^n$ be a compact region with a smooth boundary
and $g_0$ the standard hyperbolic metric on $D$.
Then there is a neighborhood $\mathcal U$ of $g_0$
in the space of Riemannian metrics on~$D$
such that for every metric $g\in\mathcal U$
the space $(D,g)$ is a minimal filling and boundary rigid.
\end{theorem}

As explained above, it suffices to prove that the metric $g$
in question is a unique minimal filling of its boundary
distance function. The space of Riemannian metrics in these
theorems is regarded with $C^\infty$ topology.
(In fact, one can lower it down to $C^2$ in Theorem \ref{theorem1}
and to $C^3$ in Theorem \ref{theorem2}.)

\subsection{Isometric representations}
It is well known that every metric space $X$ can
be isometrically embedded into an $L^\infty$ type
Banach space. A classic Kuratowski map embeds
a bounded metric space $X$ into $C^0(X)$ by
sending every point $x\in X$ to the distance function
$d_X(x,\cdot)\in C^0(X)$.
For simple Riemannian metrics there
are other natural constructions.

Let $M$ be a simple Riemannian manifold and $S=\partial M$.
The \textit{boundary distance representation}
is a map $\Phi:M\to C^0(S)\subset L^\infty(S)$
defined by
$$
 \Phi(x) = d_M(x,\cdot)|_{S} .
$$
It is easy to see that this map is distance-preserving.
Furthermore, it features additional nice properties:
it is smooth away from the boundary and the gradients
of its ``coordinate functions'' $d_M(\cdot,s)$, $s\in S$,
at every point $x\in M\setminus\partial M$ define
a diffeomorphism between $S$ and the unit tangent bundle at~$x$.
This technical property plays an important role.

There is a similar construction for
a complete simply connected manifold $M$ of nonpositive
curvature (or a compact region in such a manifold).
Fix a point $o\in M$ and let $S=UT_oM$ be the unit
tangent sphere at $o$. The \textit{Busemann representation}
$\Phi:M\to L^\infty(S)$ is defined by
\begin{equation}
\label{e-busemann}
 \Phi(x)(v) = B_{\gamma_v}(x), \qquad x\in M,\ v\in S,
\end{equation}
where $\gamma_v$ is the geodesic ray from $o$ defined by
the initial data $\dot\gamma_v(0)=v$, and $B_{\gamma_v}$
is its Busemann function.
In the case $M=\R^n$ this map is linear:
$$
 \Phi(x) = \langle x,\cdot\rangle|_{S^{n-1}}, \qquad x\in\R^n,
$$
where $S^{n-1}$ is the standard unit sphere in $\R^n$.
It is easy to see that the Busemann representation of
a nonpositively curved metric is distance-preserving.
If the metric has constant curvature outside a compact set,
then the Busemann representation is smooth
(in general, it may fail to be smooth even in the co-compact case).

The proofs of the above theorems are based on the following fact:

\begin{theorem}[\cite{I09}]
\label{theorem-aux}
Let $M$ be a compact Riemannian manifold with boundary,
$S$ a $\sigma$-finite measure space and $\Phi:M\to L^\infty(S)$ a distance-preserving map.
Then $M$ is a minimal filling if and only if
$\Phi(M)$ is an area minimizer, that is, it has the least area among all
Lipschitz surfaces in $L^\infty(S)$ with the same boundary.

Furthermore, if $\Phi(M)$ is a unique area minimizer spanning its boundary,
then $M$ is a unique minimal filling of its boundary distance function
and hence is boundary rigid.

Here the surface area in $L^\infty(S)$ is defined as the Loewner area, see below.
\end{theorem}

The proof of Theorem \ref{theorem-aux} is similar to the argument in
section \ref{sec-filling} showing that Conjecture \ref{c-main}
is equivalent to Question \ref{q-basic}.
The ``if'' implication and the uniqueness assertion easily
follow from the fact that any filling $M'$
of $(\partial M,bd_M)$ admits a 1-Lipschitz map
$\Phi':M'\to L^\infty(S)$ such that $\Phi'|_{\partial M}=\Phi|_{\partial M}$.
This part of the proof works for any definition
of surface area satisfying the natural requirement
that 1-Lipschitz maps do not increase areas.

The ``only if'' implication is not used in theorems \ref{theorem1}
and \ref{theorem2} but it is important for motivation.
This implication requires a careful choice of the surface area definition,
see the next subsection.

\begin{remark}
Theorem \ref{theorem-aux} is a partial case of the following fact.
Let $N$ be a closed  $(n-1)$-manifold,
$d:N\times N\to\R$ is a metric on~$N$ and $\Psi$
a distance-preserving map from $(N,d)$ to $L^\infty(S)$.
Then $\FillVol(N,d)$ equals the filling area of $\Psi(N)$
in $L^\infty(S)$, i.e.\ the infimum of the (Loewner) areas of Lipschitz
$n$-surfaces in $L^\infty(S)$ whose boundaries
are parametrized by~$\Psi$.

In his founding paper \cite{gromov-frm} Gromov used
the fact that filling volumes and filling areas in $L^\infty$ are equal
up to a factor bounded by a constant depending on~$n$.
This factor could not be removed because Gromov used
another definition of area
(namely Benson's area, cf.\ \cite{thompson} and \cite{benson},
denoted by $mass*$ in \cite{gromov-frm}). 
If one is interested in filling volumes up to a bounded factor,
any definition of area works fine, and $mass*$ is 
technically easier than other definitions.
However it is not suitable for finding precise
filling volumes.
\end{remark}

\subsection{Defining the surface area in $L^\infty$}
There are two issues to sort out. First, we have to deal
with surfaces of only Lipschitz regularity.
For Lipschitz surfaces in $\R^n$ one uses
Rademacher's theorem asserting that
every Lipschitz map is differentiable almost everywhere.
This gives one a Jacobian defined a.e.\ 
and then the surface area is defined by integration.
This scheme does not work for surfaces
in $L^\infty$ due to the lack of Rademacher's theorem.
This can be worked around by using \textit{weak derivatives}
(i.e., derivatives with respect to a weak topology on the target space).
For a Lipschitz map from a smooth manifold $M$ to $L^\infty$,
weak derivatives exist and have natural metric properties almost everywhere on~$M$,
cf.\ \cite{AK00} or~\cite{I09}.
(This Rademacher-type theorem is the main reason why we prefer $L^\infty$
over $C^0$ for the target space of our embeddings.)
Then, in order to define the surface area in $L^\infty$,
one uses weak derivatives in the same way as ordinary derivatives
in~$\R^n$.

The second issue is how to define the area integrand. Since the norm in $L^\infty$
is not Euclidean, the induced metric of a surface (even of a smooth one)
is not Riemannian in general. In fact, it can be an arbitrary Finsler metric.
Contrary to the Riemannian case, there are many non-equivalent definitions
of area and volume for Finsler metrics, see e.g.~\cite{thompson}.
The most commonly used definitions are Busemann's \cite{busemann47}
(the Hausdorff measure) and Holmes--Thompson's \cite{holmes-thompson}
(the projection of the Liouville measure from the unit tangent bundle).

In order to define an $n$-dimensional Finsler volume,
one chooses a volume normalization factor in every
(affine type of) $n$-dimensional Banach space.
For example, Busemann's definition normalizes
the volume of the norm's unit ball to be the same constant
$\omega_n$ for all $n$-dimensional Banach spaces.
The \textit{Loewner volume} mentioned in Theorem \ref{theorem-aux}
is defined as follows. Let $(V,\|\cdot\|)$ be an
$n$-dimensional Banach space,
$B$ its unit ball and $E$ the John--Loewner ellipsoid of $B$
(i.e., the ellipsoid of maximal volume contained in~$B$).
Then the Loewner volume in $(V,\|\cdot\|)$ is normalized so that the volume
of $E$ equals~$\omega_n$.
For a Finsler manifold $M=(M,\varphi)$, the Loewner volume equals
the infimum of volumes of Riemannian metrics $g$ on $M$
satisfying $g(v,v)\ge \varphi^2(v)$ for all $v\in TM$.
This definition extends to Lipschitz surfaces in $L^\infty$
as explained above.

\begin{remark}
Theorem \ref{theorem-aux} is valid in a
more general context of Finslerian minimal fillings.
To define the notion of a Finslerian minimal filling,
modify Definition \ref{d-filvol} of filling volume
so that the infimum in \eqref{e-fillvol}
is taken over Finsler manifolds $M$ rather than Riemannian ones.
Naturally one has to choose a definition of Finsler volume in \eqref{e-fillvol},
and the same definition should be used for the surface area
in Theorem \ref{theorem-aux}.
Choosing Loewner's volume definition yields the
Riemannian version of the theorem as a special case
of the Finslerian one, cf.~\cite{I09}.
\end{remark}

\subsection{Sketch-proof of theorems \ref{theorem1} and \ref{theorem2}}
First I explain how the proof works in the (well-known) case
when $g=g_0$, that is, $M$ is a compact region $D\subset\R^n$
equipped with the Euclidean metric.

Let $S=S^{n-1}$ and $\Phi_0:\R^n\to L^\infty(S)$ be the Busemann representation
of the standard Euclidean metric. That is, $\Phi_0$ is a linear map
defined by
\begin{equation}
\label{e-linear-embedding}
 \Phi_0(x) = \langle x,\cdot\rangle|_{S}, \qquad x\in\R^n,
\end{equation}
where $S$ is identified with the unit sphere in~$\R^n$.
Denote $W=\Phi_0(\R^n)$ and $B=\Phi_0(D)$.
By Theorem \ref{theorem-aux} it suffices to prove that $B$
is a unique Loewner area minimizer
in $L^\infty(S)$ among the Lipschitz surfaces with the same boundary.
In fact, we can restrict ourselves to surfaces contained in
a sufficiently large ball.

Equip $L^\infty(S)$ with a scalar product $\langle\cdot,\cdot\rangle_e$ defined by
by
\begin{equation}
\label{e-scalar-product}
 \langle u,v\rangle_e = n\int_S uv\,d\mu
\end{equation}
where $\mu$ is the Haar probability measure on~$S$. 
This defines a Euclidean norm on $L^\infty(S)$
that we denote by $\|\cdot\|_e$.
One easily sees that $\|\cdot\|_e$ is Lipschitz
w.r.t.\ the $L^\infty$ norm and the two norms coincide on~$W$.
An easy application of Cauchy--Schwartz
inequality shows that the Euclidean $n$-volume defined by the above
scalar product is no greater than the Loewner $n$-volume defined
by the $L^\infty$ norm.
Hence the Euclidean $n$-area of any Lipschitz surface
in $L^\infty(S)$ is no greater that the Loewner $n$-area,
and these areas are equal if the surface is contained in~$W$.
Thus it suffices to prove that $\Phi_0(D)$ minimize the Euclidean
area among the surfaces with the same boundary.
And this is trivial because the orthogonal projection
onto $W$ (with respect to our scalar product) does not increase areas.

Furthermore, one can compose the projection with a suitable shrinking in $W$
to obtain a smooth retraction $P:L^\infty(S)\subset L^2(S)\to W$ 
such that
\begin{equation}
\label{e-jacobian}
 J_nP(u) \le 1-c\cdot\|u-P(u)\|_e^2
\end{equation}
for some $c>0$ and all $u$ from a large ball in $L^2(S)$.
Here $J_n$ denotes the $n$-dimensional Jacobian with respect
to $\|\cdot\|_e$.
This proves uniqueness and a sort of stability estimate.

Now consider the general case of Theorem \ref{theorem1}
when the metric $g$ of $M=(D,g)$ is close
to Euclidean in $C^r$ topology for a suitable $r$ (in fact, $r=3$
is sufficient for the argument presented here and
a more delicate argument in \cite{BI10} works for $r=2$).
The proof of Theorem \ref{theorem1} consists of three steps.

\medbreak\textit{Step 1}. Construct a smooth distance-preserving
map $\Phi:M\to L^\infty(S)$ close to the above linear map $\Phi_0$
(in a suitable topology).
In order to do this, one can use a formula similar to \eqref{e-busemann}
with Riemannian distances to hyperplanes rather than Busemann functions.
By Theorem \ref{theorem-aux}, it suffices to prove that $\Phi(M)$
is a unique Loewner area minimizer among the surfaces with the same boundary.

\medbreak\textit{Step 2}. Prove that the surface $\Phi(M)$ is minimal
in a variational sense. This part of the proof is the most encouraging:
it does not depend on the fact that the metric
is close to Euclidean and works for any boundary distance representation
of a simple metric, any smooth Busemann representation and, in fact,
for any isometric embedding with a similar behavior of coordinate functions.

What is meant by being a minimal surface needs clarification.
Unfortunately, the first variation of the Loewner area does not
make sense since the Loewner area integrand is not differentiable
(even in a finite-dimensional Banach space with a smooth norm).
To work around this, we differentiate a smooth lower bound for the Loewner area.
This lower bound is the $n$-area defined by a Riemannian metric $\mathcal G$
on $L^\infty(S)$ extending the metric of $\Phi(M)$. 

The metric $\mathcal G$ is a smooth family of scalar products
$\langle\cdot,\cdot\rangle_\varphi$, $\varphi\in L^\infty(S)$, on $L^\infty(S)$.
Every scalar product $\langle\cdot,\cdot\rangle_\varphi$ is given
by a formula similar to \eqref{e-scalar-product}
where $\mu$ is replaced by a probability measure $\mu_\varphi$ depending on $\varphi$.
The normalization of the measures $\mu_\varphi$ implies that
the $n$-area defined by $\mathcal G$ is no greater than
the Loewner $n$-area.
In order to make $\mathcal G$ compatible with the metric
of $\Phi(M)$, one defines $\mu_\varphi$
explicitly for every $\varphi\in\Phi(M)$.
Namely if $\varphi=\Phi(x)$ where $x\in M$, then the measure $\mu_\varphi$
is obtained from the normalized Haar measure on the unit sphere $UT_xM\subset T_xM$
via a natural diffeomorphism between $UT_xM$ and~$S$.
(This diffeomorphism turns the derivative $d_x\Phi:T_xM\to L^\infty(S)$
into the standard linear map given by \eqref{e-linear-embedding}).

The variational minimality of $\Phi(M)$ means that the first variation
of the Riemannian $n$-area defined by $\mathcal G$
is zero for every (Lipschitz) variation, or, equivalently,
the mean curvature w.r.t.\ any normal vector is zero.
The proof is a direct computation of the mean curvature.
It works for any Riemannian structure $\mathcal G$ defined as above,
however the next step assumes that $\mathcal G$ is
a small perturbation of the flat Riemannian structure
defined by~\eqref{e-scalar-product}.

\medbreak\textit{Step 3}. Prove that $\Phi(M)$ is a unique
area minimizer with respect to~$\mathcal G$ provided  that
$\mathcal G$ is sufficiently close (in a suitable topology)
to the constant scalar product \eqref{e-scalar-product}.
Since the $n$-area defined by $\mathcal G$ is a lower bound for
the Loewner $n$-area and the two areas coincide on $\Phi(M)$,
it follows that $\Phi(M)$ is a unique minimizer of the Loewner
area and hence $M$ is a minimal filling and boundary rigid.

The proof essentially establishes the fact that stable minimality
that we had in the case $g=g_0$ is stable under small perturbations
of the data. More precisely, one can construct a retraction
from $L^\infty(S)$ to a (minimal) surface containing $\Phi(M)$
by perturbing the area-decreasing map $P$ used in the flat case.
The perturbation should preserve the property that pre-images of
points are orthogonal to the surface. Since the surface is minimal,
this implies that the $n$-dimensional Jacobian
(with respect to $\mathcal G$) of the
retraction has zero derivatives at the surface.
And if its second derivatives are close to the original
ones, the inequality \eqref{e-jacobian} persists,
implying the desired result.

\medbreak\textit{Proof of Theorem \ref{theorem2}}.
The proof goes along the same lines: first we prove
the desired properties for the standard hyperbolic metric
and then verify that they are stable under perturbations.

The only essential difference is the choice of an area
non-increasing map in place of the linear orthogonal projection.
We define a ``projection'' $P:L^\infty(S)\to\HH^n$ as follows:
for every $\varphi\in L^\infty(S^{n-1})$,
$P(\varphi)$ is a (unique)
point where the function $F_\varphi:\HH^n\to\R$ defined by
$$
F_\varphi(x)=\int_S e^{-n\varphi(s)}e^{B_{\gamma_s}(x)} \,ds
$$
attains its minimum.
Here $S=T_o\HH^n$ where $o\in\HH^n$ is a fixed origin,
$B_{\gamma_s}$ denotes the Busemann function of a geodesic ray
starting from the origin in the direction~$s$,
and $ds$ denotes the standard measure on~$S$.

One can verify that $P$ does not increase $n$-dimensional Loewner areas
and that $\Phi_0\circ P$ is a retraction of $L^\infty(S)$
onto $\Phi_0(\HH^n)$ where $\Phi_0$ is the Busemann representation
of~$\HH^n$. This proves filling minimality and boundary rigidity
for regions in~$\HH^n$.
Then the proof of Theorem \ref{theorem2} is similar to that
of Theorem \ref{theorem1}.

\section{Finslerian case}

As shown by Theorem \ref{theorem-aux}, reducing
filling minimality to area minimality is a natural
approach (at least there is no loss of generality
at this step).
But some other tricks in the above proofs are too limited;
it would be nice to replace them by
a better technique.
In particular, replacing the Loewner area by the area
defined by an auxiliary Riemannian metric $\mathcal G$
is suspicious: this may not work for other minimal
fillings, and there is no natural way to choose
this auxiliary metric.

It would be more natural to utilize the Finslerian
nature of surfaces in $L^\infty$ and work with
their natural Finsler areas, e.g.\ Busemann
or Holmes--Thompson areas. Unfortunately very little
is known about these surface areas in co-dimensions higher than~1.
For example, the following basic question is not yet
answered.

\begin{question}[Busemann \cite{busemann61}]
\label{q-elliptic}
Let $V$ be a finite-dimensional Banach space,
$D$ an $n$-disc in an $n$-dimensional affine
subspace $W\subset V$ and $F$ is an orientable surface
in $V$ such that $\partial F=\partial D$.
Is it always true that $\area(F)\ge\area(D)$?
\end{question}

In other words, is the $n$-dimensional area integrand
in a Banach space semi-elliptic (over $\Z$)?
Actually this is a different question
for every definition of area.
In the cited paper \cite{busemann61} the
question is asked for the Holmes--Thompson area,
defined there in terms of the projection function
of a convex body.
For both Busemann and Holmes--Thompson areas,
the answer is known to be affirmative in the case $\dim V=n+1$
but the question is open in higher co-dimensions
(even in the special case when the restriction of the Banach norm
to the subspace~$W$ is Euclidean). Contrary to this,
Benson area and Loewner area are known to
be semi-elliptic in all dimensions and co-dimensions,
cf.\ \cite{gromov-frm} and \cite{I09}.

An affirmative answer to Question \ref{q-elliptic} would
have nice applications including a Finslerian generalization
of the asymptotic volume estimate \eqref{e-asvol},
cf.~\cite{BI02}.
It would also imply that every region in an $n$-dimensional
Banach space is a Finslerian minimal filling.
This is especially interesting in the case of the Busemann
volume because it is equal to the Hausdorff measure
naturally defined for all metric spaces, not just Finslerian.
Here is how one can formulate a filling question without
referencing anything from differential geometry.

\begin{question}
\label{q-hausdorff}
Let $d$ be a (continuous) metric on the standard unit ball
$D^n\subset\R^n$ such that
$$
d(x,y)\ge d_E(x,y):=|x-y|
$$
for all $x,y\in\partial D^n=S^{n-1}$. Is it true that
for all such metrics $d$ one has
$$
\mathcal H^n(D^n,d) \ge \mathcal H^n(D^n,d_E)
$$
where $\mathcal H^n$ denotes
the $n$-dimensional Hausdorff measure?
\end{question}

An affirmative answer to Question \ref{q-elliptic}
would answer Question \ref{q-hausdorff} for a Lipschitz
metric~$d$. I do not know whether the case of a
general metric is different.

One may also seek a Finslerian generalization
of the minimal filling conjecture \ref{c-main}.
Although there is no boundary rigidity in the Finslerian case,
simple Finsler metrics sharing the same boundary distance function
have the same Holmes--Thompson volume. 
This leaves a possibility that the Finslerian generalization
of Conjecture \ref{c-main}
might be true if the volume of a Finsler
metric is defined as the Holmes--Thompson volume.
This generalization is ``almost proved'' in dimension~2:
every simple Finsler metric on $D^2$ is a minimal filling
\textit{among the Finsler fillings homeomorphic to $D^2$},
cf.\ \cite{I01} and~\cite{I09-pre}.

This implies a partial answer to Question \ref{q-elliptic}
for $n=2$: an affine 2-disc in a Banach space
minimizes the Holmes--Thompson area
among the surfaces spanning the same boundary and
homeomorphic to~$D^2$.
On the other hand, one can construct a Banach norm in $\R^4$
such that the resulting two-dimensional Holmes--Thompson
area integrand is not convex (that is, it
does not admit a convex extension to
the exterior product $\Lambda^2\R^4$),
cf.\ \cite{BES1}, \cite{BI02}.
And this implies that there is an affine 2-disc
which does {\em not} minimize the Holmes--Thompson area
among the Lipschitz (or polyhedral)
chains with rational coefficients, cf.~\cite{BI04}.
What is not known is whether an affine 2-disc
minimizes area among the chains with integer coefficients,
or, equivalently, among the orientable surfaces
of arbitrary genus.

\end{document}